\documentclass{amsart}

\usepackage{graphicx}
\usepackage{charter}
\usepackage{latexsym,amsfonts} 
\usepackage{amsmath,amssymb,graphics,setspace}
\usepackage{amsthm}
\usepackage{mathrsfs}
\usepackage{paralist}

\usepackage{color}
\usepackage{marvosym}

\newtheorem{example}{Example}
\newtheorem{remark}{Remark}

\makeatletter
\@namedef{subjclassname@2020}{\textup{2020} Mathematics Subject Classification}
\makeatother

\begin{document}
\title[Framework for Fractional Stochastic Equations]{
Framework for Solving Fractional Stochastic\\ 
Integral-Differential Equations}	

\author[Omid. Taghipour Birgani]{Omid T. Birgani$^{\alpha*}$}
\address{\llap{$^{\alpha}$}
Department of Mathematics, 
	University of Science and Technology of Iran, Hengam Street, Resalat Square, Tehran, Iran.//(Email Address:o.taghipour69@gmail.com)}
	
\author[J.F. Peters]{James F. Peters$^{\beta}$}
\address{\llap{$^{\beta}$}
	Department of Electrical \& Computer Engineering,
	University of Manitoba, WPG, MB, R3T 5V6, Canada and
	Department of Mathematics, Faculty of Arts and Sciences, And\.{i}yaman University, 02040 Ad\.{i}yaman, Turkey.//(Email Address: James.Peters3@umanitoba.ca)}
	
\author{Sareh Kouhkani$^{\gamma}$}
\address{\llap{$^{\gamma}$}\,
Department of Mathematics,
	Islamic Azad University branch of Shabestar,East Azerbaijan, Shabestar, Iran.//(Email Address: skouhkani@yahoo.com)}

\thanks{*Corresponding author.}



\begin{abstract}
This article introduces a framework for measuring the uncertain behaviour of a changing system in terms of the solution of a class of fractional stochastic differential equations (fsDEs). This is accomplished via operational matrices based on 2-dimensional shifted Legendre polynomials. By using operational matrices, an fsDE is converted into a matrix form and the numerical solution of the represented motion system is then found.\\
\end{abstract}


\maketitle
\tableofcontents

\section{Introduction}
\begin{figure}
	\centering
	\includegraphics[width=85mm]{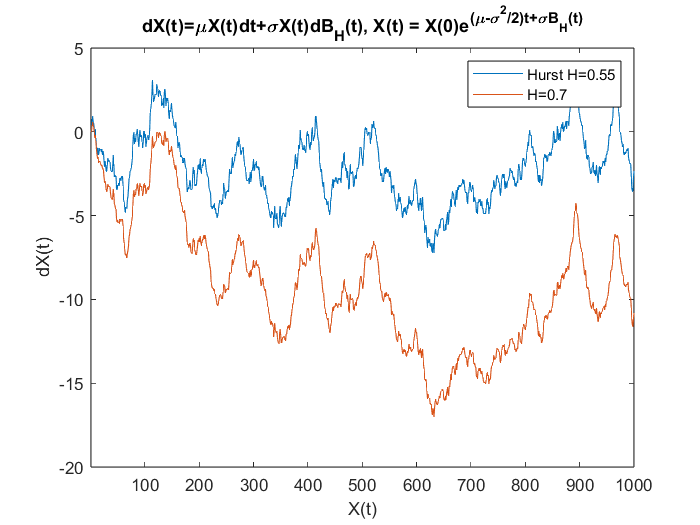}
	\label{fig:sde}
	\caption{sDE for Brownian Motion Sampling}
\end{figure}

This paper gives a numerical solution for the fractional stochastic integro-differential function
$D^\alpha f(t)$

\begin{center}
\vspace*{0.2cm}
\boxed{\boldsymbol{
D^\alpha f(t)=g(t)+ \int_{0}^{t} k_1(s,t) f(s) ds+ \int_{0}^{t} k_2(s,t) f(s) dB(s).
}}\label{fig:fSDE}
\end{center}
\vspace*{0.2cm} 

\noindent for functions $g(t)$, $k_1(t,s)$ and $k_2(t,s)$ on the probability space $(\Omega,F,P)$ and function $f(t)$ is an unknown function, $B(t)$ is a Brownian motion such that $t,s \in [0,1]$.  This is accomplished by using operational matrices based on the 2D-shifted Legendre polynomials, which are widely used in many areas of numerical analysis[33-18].

\begin{remark} 
Two sample stochastic differential equation (sDE) waveforms resulting from the Brownian motion $B_H(t)$ represented in the sDE $dX(t)$ are shown in Fig.~\ref{fig:fSDE}. The Hurst parameter $H$ used in $B_H(t)$ measures the long term memory of the time series $X(t)$[45] and is represented by the horizontal axis in Fig.~\ref{fig:fSDE}.  $X(t)$ is defined by

\begin{center}
\vspace*{0.2cm}
\boxed{\boldsymbol{
X(t)= X(0)e^{(\frac{\mu-\sigma^2}{2})t+\sigma B_H(t)}.
}}
\end{center}
\vspace*{0.2cm} 

For the derivation of the $X(t)$ in Fig.~\ref{fig:fSDE}, see [46].
\qquad \textcolor{blue}{$\blacksquare$}
\end{remark}

Univariate and bivariate orthogonal polynomials provide an algorithmic approach to convert many differential problems to a system of algebraic equations. In recent years, multivariate orthogonal polynomials have received attention in dealing with various problems have been used for solving numerical partial differential equations by many authors. In [14], I. Hesameddini and M.Shahbazi proposed the 2-dimensional shifted Legendre polynomials for solving the fractional integral equations. In [35], Shehili et al explained the bivariate Legendre approximation and applied the proposed method to some bivariate functions. Moreover, E.H. Doha[8] represents a double Chebyshev spectral method for the solution of Poisson's equation. In addition, G. Yuksel and M. Sezer[39] have applied a method based on bivariate Chebyshev series to approximate the linear second order PDEs with complicated conditions. In [7], Dascioglu has also presented a method based on approximation by truncated double Chebyshev equations with variable coefficients under most general conditions.

In recent years, the study of numerical solutions of stochastic differential equations[23] as well as fractional differential equations [42-44]
have been investigated. In [32], some discrete approximation methods are explained for solving stochastic differential equations. Misawa[29] proposed the formulation of composition methods for solving stochastic differential equations. In [25], Li et al proposed explicit numerical approximations for stochastic differential equations in finite and infinite horizons. Recently, efficient methods based on Runge- Kutta methods were proposed [19].

Fractional calculus is applied in various phenomena in viscoelasity, fluid mechanic, biology, control theory and other areas of science. Mathematics folklore sets the birth of the concept of fractional calculus in the year 1695 by the answer to a question raised by LoHpital(1661-1704) to Leibniz(1646-1716).
\section{Preliminaries}
In this section, we give the basic definitions, notations and fundamental assumptions that provide a foundation for the results given in this paper.
There are several approaches to the generalization of the notation of differentiation to fractional orders e.g. the Riemann-Liouville and the Caputo approach.\\
\textbf{Definition 1}. A fractional integral operator of order $\alpha(\alpha>0)$ of type Riemann-Liouville [2-24] is defined by
\vspace*{0.2cm}
\begin{equation}
\begin{cases}
J^{\alpha} f(t)=\frac{1}{\Gamma(\alpha)} \int_{0}^{t} (t-s)^{\alpha-1} f(s) ds,&\alpha,t>0,\\
J^0 f(t)=f(t).
\end{cases}
\end{equation}
\vspace*{0.2cm}
In (1), $\Gamma$ is the Gamma function. The fundamental properties of $J^\alpha$ for $f(t)$ are
\begin{compactenum}[1$^o$]
\item $J^\alpha J^\theta f(t)=J^{\alpha+\theta} f(t)$.
\item $J^\alpha t^\gamma =\frac{\Gamma(\gamma+1)}{\Gamma(\alpha+\gamma+1)} t^{\alpha+\gamma}, \Gamma(\alpha+\gamma+1) > 0$.
\item $J^\theta J^\alpha f(t)=J^\alpha J^\theta f(t)$, commutativity property.
\end{compactenum} 
\vspace*{0.2cm}
\textbf{Definition 2.} Fractional derivatives of order $\alpha$ of $f(t)$ of type Caputo are defined by
\begin{equation}
D^\alpha f(t)=J^{n-\alpha} D^n f(t)=\frac{1}{\Gamma(n-\alpha)} \int_{0}^{t} (t-s)^{n-\alpha-1} \frac{d^n}{ds^n} f(s), t>0, n-1<\alpha<n.
\end{equation}
In (2), $D^n$ is a differential operator of order $n$. In addition,
\begin{equation}
D^\alpha t^\mu =\begin{cases}
0&\mu\in N_0, \mu<\lceil \alpha \rceil \\
\frac{\Gamma(\mu+1)}{\Gamma(\mu+1-\alpha)} t^{\mu- \alpha}&\mu \in N_0, \mu\geq \lceil \alpha\rceil or \mu \notin N, \mu>\lfloor \alpha\rfloor,
\end{cases}
\end{equation}
In (3), $N=1,2,\dots$ (natural numbers) and $N_0=0,1,2,\dots$ (an integer greater than or equal to 0). Next, consider\\
\textbf{Definition 3}. Let $(\Omega,F,P)$ be a probability space with filteration ${\{F_t\}}_{t\geq 0}$. The Brownian motion $B$ is a stochastic process over a a probability space $(\Omega,F,P)$ with filteration ${\{F_t\}}_{t\geq 0}$, having the following properties [27]:
\begin{enumerate}
	\item $B(0)=0$ (0 is a Spanier-Brouwer fixed point for $B:\mathbb{R}\rightarrow \mathbb{R}$~[31],[41],[42]).
	\item $B(t)-B(s), t>s$ is independent of $F_s$, such that
	\item $B(t)-B(s)$ is normally disributed with mean $\mu = 0$ and variance $t-s$.
	\item $B(t), t\geq 0$ is a continuous function.
\end{enumerate}
\vspace*{0.2cm}
\textbf{Definition 5}. An integral of the form$\int_{a}^{b} f(t) dB(t)$ where $B(t)$ denotes the Brownian motion process can be defined in various ways. This integral can be approximated as
\begin{equation}
\sum_{i=1}^{n} f(\tau_i)(B(t_i)-B(t_{i-1}))
\end{equation}
where $\tau_i$ is chosen arbitrary in the subinterval $[t_{i-1},t_i]$ of the partition $ a=t_0<t_1<\cdots<t_n=b$. If $\tau_i$ is taken as $t_{i-1}$, (4) is an It\^o integral and if , (4) is a Stratonovich integral.
\section{Properties of Legendre and shifted Legendre polynomials}
The Legendre polynomials $L_n(x)$, are defined on the interval $I=[-1,1]$ have the following properties
\begin{equation}
L_0=1,
 L_1=x,
 \cdots,
  L_\text{i+1}(x)=\frac{2i}{i+1}L_i(x)-\frac{i}{i+1}L_\text{i-1}(x),
\end{equation} 
In order to use these polynomials on the interval $x\in I=[a,b]$ we define the so-called shifted Legendre polynomials by introducing the change of variable $t=\frac{2x}{I}-1$. Let $P_i(x)$ denote the shifted Legendre polynomials and we define shifted Legendre polynomials as:
\begin{equation}
\begin{cases}
P_0=1, P_1(x)=\frac{2x-(b+a)}{b-a}\\
P_i(x)=(2-\frac{1}{i})(\frac{2x-(b+a)}{b-a})P_\text{i-1}(x)-(1-\frac{1}{i})P_\text{i-2}(x), i=2\cdots,
\end{cases}
\end{equation}
For example these polynomials on the interval $[0,1]$ can be obtained as 
\begin{equation}
P_\text{i+1}(x)=\frac{(2i+1)(2x-1)}{(i+1)}P_i(x)-\frac{i}{i+1}P_\text{i-1}(x), i=1,2,\cdots,
\end{equation}
where $P_0=1$ and $P_1(x)=2x-1$.
The analytical form of the shifted Legendre polynomials $P_i(x)$ of degree $i$ is given by
\begin{equation}
P_i(x)=\sum_{k=0}^i (-1)^\text{i+k} \frac{(i+k)!x^k}{(i-k)!(k!)^2 I^k}
\end{equation}
note that $P_{i}(0)=(-1)^i$ and $P_i(1)=1$. The orthogonality condition is
\begin{equation}
\int_0^I P_i(x)P_j(x)\omega_I(x)dx=h_I,k\sigma_jk,
\end{equation}
where $\omega_I=1$ and $h_I,k=\frac{I}{(2k+1)}$.
A function $u(x)$, square integrable on $I=[a,b]$, can be expressed in terms of shifted Legendre polynomials as
\begin{equation}
u(x)=\sum_{j=0}^{m} u(x)P_j(x)dx, j=0,1,2,\cdots.
\end{equation}
In practice, only the first $(m+1)$-terms of shifted Legendre polynomials are considered. Hence, we have
\begin{equation}
u(x)=\sum_{j=0}^{m} c_jP_(x)=C \psi(x),
\end{equation}
where $C$ is a shifted Legendre coefficient vector and $\phi$ is a shifted Legendre vector as the following
\begin{equation}
C=[c_0,c_1,\cdots,c_m],  \psi(x)=[P_0(x),P_1(x),\cdots,P_m(x)]^T.
\end{equation}
\textbf{Theorem 1}. The derivative of the vector $\psi(x)$ can be expressed by
\begin{equation}
D\psi(x)=\Lambda^\text{(1)} \psi(x),
\end{equation}
where $\Lambda^{(1)}$ is the $(m+1) \times(m+1)$ operational matrix of derivative given by
\begin{equation}
\Lambda^{(1)}=[(d_{ij})]^T=
\begin{cases}
$2(2j+1), \text{for} j=i-k$,\\
0, \text{otherwise},
\end{cases}
\end{equation}
where $i,j$ are numbers of the rows and columns respectively on the interval $[0,1]$.\\
For example for even $m$ we have
\begin{equation*}
\Lambda^{(1)}=\begin{bmatrix}
0&0&0&0&\ldots&0&0&0\\
1&0&0&0&\ldots&0&0&0\\
0&3&0&0&\ldots&0&0&0\\
1&0&5&0&\ldots&0&0&0\\
1&0&5&0&\ldots&(2m-3)&0&0\\
0&3&0&7&\ldots&0&(2m-1)&0
\end{bmatrix}
\end{equation*}
By using equation(8) for $n \in N$, we have
\begin{equation}
D^n \psi(x)=(\Lambda^{(1)})^n \psi(x).
\end{equation}\\
\textbf{Theorem 2.} Let $\psi(x)$ be a shifted Legendre vector with $\alpha>0$, then, from [33], we have\\
$D^\alpha \psi(x)\simeq D^{(\alpha)} \psi(x)$,\\
such that $D^{(\alpha)}=(m+1)\times (m+1)$, an operational matrix of a Caputo fractional derivative of order $\alpha$ defined by
\begin{equation}
D^{(\alpha)}=\begin{bmatrix}
0&0&\ldots&0\\
\vdots&\vdots&\ldots&\vdots\\
0&0&\ldots&0\\
\sum_{k=\lceil \alpha \rceil}^{\lceil \alpha \rceil} \theta_{\lceil \alpha \rceil,0,k}&\sum_{k=\lceil \alpha \rceil}^{\lceil \alpha \rceil} \theta_{\lceil \alpha \rceil,1,k}&\ldots&\sum_{k=\lceil \alpha \rceil}^{\lceil \alpha \rceil} \theta_{\lceil \alpha \rceil,m,k}\\
\vdots&\vdots&\ldots&\vdots\\
\sum_{k=\lceil \alpha \rceil}^{i} \theta_{i,0,k}&\sum_{k=\lceil \alpha \rceil}^{i} \theta_{i,1,k}&\ldots&\sum_{k=\lceil \alpha \rceil}^{i} \theta_{i,m,k}\\
\vdots&\vdots&\ldots&\vdots\\
\sum_{k=\lceil \alpha \rceil}^{m} \theta_{m,0,k}&\sum_{k=\lceil \alpha \rceil}^{m} \theta_{m,1,k}&\ldots&\sum_{k=\lceil \alpha \rceil}^{m} \theta_{m,m,k}\\
\end{bmatrix},
\end{equation}
where 
\begin{equation}
\theta_{i,j,k}=(2j+1)\sum_{l=0}^{j} \frac{(-1)^{i+j+k+l} (i+k)! (l+j)!}{(i-k)! k! \Gamma{(k-\alpha+1)} (j-l)! (l!)^2 (k+l-\alpha +1)}.
\end{equation}\\
\textbf{Proof}. See[33].
\\2-D shifted Legendre polynomials are defined on $\Delta=[0,a] \times[0,b]$ as follows
\begin{equation}
P_{r,s}(x,t)=P_r(x) P_s(t)=\psi(x) \times \Phi(t) , i,j=0,1,2,\cdots,
\end{equation}
and from  [30], $P_{r,s}(x,t)$ is defined by
\begin{equation}
\int_{0}^{s} \int_{0}^{z} P_{r,s}(x,t) dxdt \simeq W P_{r,s}(x,t)=(Q_{1} \otimes Q_{2}) P_{r,s}(x,t)
\end{equation}
such that $x \in[0,a], t \in[0,b]$, $W$ is $(m+1)^2 \times (m+1)^2$
and $Q_1$, $Q_2$ are operational matrices of 1-D shifted Legendre polynomials on $[0,a]$ and $[0,b]$ defined by
\begin{align*}
Q_1 &=\begin{bmatrix}
1&1&0&\ldots&0&0&0\\
-\frac{1}{3}&0&\frac{1}{3}&\ldots&0&0&0\\
\vdots&\vdots&\vdots&\vdots&\vdots&\vdots&\vdots\\
0&0&0&\ldots&-\frac{1}{2m-1}&0&\frac{1}{2m-1}
\end{bmatrix},\\
Q_2 &=\begin{bmatrix}
1&1&0&\ldots&0&0&0\\
-\frac{1}{3}&0&\frac{1}{3}\ldots&0&0&0\\
\vdots&\vdots&\vdots&\vdots&\vdots&\vdots&\vdots\\
0&0&0&\ldots&-\frac{1}{2m-1}&0&\frac{1}{2m-1}\\
0&0&0&\ldots&0&-\frac{1}{2m+1}&0
\end{bmatrix}.
\end{align*}
Also, we can write[14]
\begin{equation}
\int_{0}^{x} P(x,t) dx \simeq Q_{3} P(x,t), 
\int_{0}^{t} P(x,t) dt \simeq Q_{4} P(x,t).
\end{equation}\\
where $Q_{3}$ and $Q_{4}$ are $(m+1)^{2} \times (m+1)^{2}$ matrices of the form
\begin{equation*}
Q_3=\begin{bmatrix}
I&I&O&\ldots&O&O&O\\
-\frac{I}{3}&O&\frac{I}{3}&\ldots&O&O&O\\
\vdots&\vdots&\vdots&\ddots&\vdots&\vdots&\vdots\\
O&O&O&\ldots&-\frac{I}{2m-1}&O&\frac{I}{2m-1}\\
O&O&O&\ldots&O&-\frac{I}{2m+1}&O
\end{bmatrix},
Q_4=\begin{bmatrix}
p_2&O&O&\ldots&O\\
O&p_2&O&\ldots&O\\
O&O&p_2&\ldots&O\\
\vdots&\vdots&\vdots&\ddots&\vdots\\
O&O&O&\ldots&p_2
\end{bmatrix}.
\end{equation*}\\
where I and O are the identify and zero matrix of order $(m+1)$, respectively and $p_2$ is a square matrix of order $(m+1)$ which is defined by
\begin{equation*}
p_2=\frac{b-a}{2} \begin{bmatrix}
1&1&0&\ldots&0&0&0\\
-\frac{1}{3}&0&\frac{1}{3}&\ldots&0&0&0\\
\vdots&\vdots&\vdots&\ddots&\vdots&\vdots&\vdots\\
\vdots&\vdots&\vdots&\ddots&-\frac{1}{2m+1}&0&\frac{1}{2m+1}\\
0&0&0&\ldots&0&-\frac{1}{2m+1}&0
\end{bmatrix}.
\end{equation*}

\textbf{Theorem 3}. We can express the function of $u(x,t)$ as the following matrix form
\begin{equation}
u(x,t)=\sum_{r=0}^{\infty} \sum_{s=0}^{\infty} c_{r,s} P_{r,s}(x,t)=\psi(x) \Phi(t) \bar{C},
\end{equation}
where $\bar{C}$ is the unknown Legendre coefficient matrix
\begin{equation*}
\bar{C}=[c_{0,0},c_{0,1},\cdots,c_{0,m},c_{1,0},c_{1,1},\cdots,c_{1,m},\cdots,c_{m,0},\cdots,c_{m,m}]_{(m+1)^2 \times1}^{T}
\end{equation*}
where the coefficient matrix $\bar{C}$ is given by
\begin{equation}
c_{i,j}=\frac{(2i+1)}{b-a} \int_{a}^{b} \int_{a}^{b} u(x,t)L_{i}(x) L_{j}(t) dx dt
\end{equation}
and
\begin{equation*}
\Phi(t)=\begin{bmatrix}
P_0(t)&\ldots&P_m(t)&0&0&0&\ldots&0&0&0\\
0&0&0&P_0(t)&\ldots&P_m(t)&\ldots&0&0&0\\
\vdots&\vdots&\vdots&\vdots&\vdots&\vdots&\ddots&\vdots&\vdots&\vdots\\
0&0&0&0&0&0&\ldots&P_0(t)&\ldots&P_m(t)\\  
\end{bmatrix}.
\end{equation*}
is a matrix of order ${(m+1)\times (m+1)^2}$.\\
\textbf{Proof}: See proof in[24].\\

\section{Analysis of the solution method}
This section is devoted to analyze a numerical solution of the following fractional stochastic integro-differential equation based on the presented method
\begin{equation}
D^\alpha f(t)=g(t)+ \int_{0}^{t} k_1(s,t) f(s) ds+ \int_{0}^{t} k_2(s,t) f(s) dB(s),
\end{equation} 
where $g(t)$, $k_1(t,s)$ and $k_2(t,s)$ are given functions on the same probability space $(\Omega,F,P), f(t)$ is the unknown function and $B(t)$ is a Brownian motion also, $t,s \in [0,1]$ are variables.\\
In order to use the proposed method, we pursue the following process. We expand functions $f(t), g(t)$ in the following approximations
\begin{equation}
f_n (t)=\bar{f}\psi(t),
\end{equation}
\begin{equation}
g_n (t)=\bar{g}\psi(t),
\end{equation}
where $\bar{f}$ and $\bar{g}$ are coefficient vectors that $ \bar{f}=[f_{0}, f_1,\cdots,f_m] $ and $ \bar{g}=[g_0, g_1, \cdots, g_m]$. Also, by using equation(21) we approximate $k(s,t)$ as the following
\begin{equation}
k_n (s,t)=\psi(s) \Phi(t) \bar{k},
\end{equation}
where $\bar{k}$ is a vector and its elements obtain similar to equation(22)
Now, we approximate $ D^\alpha f(t)$ as
\begin{equation}
D^\alpha f_n(t)=\bar{f} D^{(\alpha)} \Psi(t). 
\end{equation}
In order to approximate the first integral in equation(23) by applying equation(20), we have
\begin{equation}
\int _{0}^{t} k_{1,m} (s,t) f_{m}(s) ds=\int_{0}^{t} \psi(s) \Phi(t) \bar{C} \Psi(s) \bar{f} ds =\bar{f} \Phi(t) Q_{3} \bar{k_1}.
\end{equation}
Now we approximate the second integral in equation(23) by applying the following operation mertix from equation (20)
\begin{equation}
\int_{0}^{t} P(s,t)dB(s)=\bar{Q_3}P(s,t),
\end{equation}
where $\bar{Q_3}$ is a matrix of order $(m+1)\times (m+1)$ is given by
\begin{equation}
\bar{Q_3}=\begin{bmatrix}
\bar{B}&\bar{B}&O&\ldots&O&O&O\\
-\frac{\bar{B}}{3}&O&\frac{\bar{B}}{3}&\ldots&O&O&O\\
\vdots&\vdots&\vdots&\ddots&\vdots&\vdots&\vdots\\
O&O&O&\ldots&-\frac{\bar{B}}{2m-1}&O&\frac{\bar{B}}{2m-1}\\
O&O&O&\ldots&O&-\frac{\bar{B}}{2m+1}&O
\end{bmatrix},
\end{equation}
where $\bar{B}$ is a diagonal matrix of order $(m+1)$ which the elements of this matrix is given by\\
$\bar{B_i}=B(\varphi_i)-B(\varphi_{i-1}) $, \\
where $B(\varphi)$ denots the Brownian motion which is normalized as follows\\
$B(\varphi)=\sqrt{\varphi} N(0,1),$
on the subinterval $[\varphi_{i-1}, \varphi_{i}]$ of the partition $a=\varphi_0<\varphi_1<\ldots<\varphi_{m}=b$.\\
Now, by employing equations(4), (26), (20) and (29) in the second integral in equation(23) we have
\begin{equation}
\int_{0}^{t} k_{2,m} (s,t) f_m(s) dB(s)=\int_{0}^{t} \psi(s) \Phi(t) \bar{k_2} \psi(s) \bar{f} dB(s)=\bar{f} \Phi(t)^{T} \bar{Q3} \bar{k_2}.
\end{equation}
where vector $\bar{k_2}$ is known and obtain similar to equation (22).\\
We collocate equation (23) at (m+1) points $s_{\delta}=t_{\delta}, \delta=0,1,2,\ldots,(m+1)$ as\\
\begin{equation}
\bar{f}D^\alpha \psi(t_{\delta})=\bar{g}\psi(t_{\delta})+ \int_{0}^{t} k_{1,m}(s_{\delta},t_{\delta}) f_m(s_{\delta}) ds+ \int_{0}^{t} k_{2,m}(s_{\delta},t_{\delta}) f_m(s) dB(s).
\end{equation}
For collocation points, we use roots of shifted Legendre polynomials $P_{(m+1)-r}(t),m>r$, these equations together with $r$ equations of boundary and initial conditions give $(m+1)$ linear or nonlinear algebraic equations. These $(m+1)$ equations can be solved by using Gauss elimination method or Newton's iterative method respectively. 
\section{Error Analysis}
In this section, an error estimation for the approximate solution is provided. For this purpose, we pursue the proposed convergence analysis in [38] by Taheri et al. \\
To obtain an error estimation for the approximate solution of equation (23) with supplementary conditions, we call $e_{n}(t)=f(t)-f_{m}(t)$ the error function of the presented method. Therefore, $f_{m}(t)$ satisfies the following equation
\begin{equation}
D^\alpha f_m(t)=g_n(t)+\int_{0}^{t} k_{1,m}(t,s)f_{m}(s) ds+\int_{0}^{t}k_{2,m}(t,s) f_{m}(s) dB(s) +H_{m}(t), t\in [0,1]
\end{equation}
which $H_{m}(t)$ is the perturbation term. Now by subtracting equation (23) from equation (30)
\begin{align*}
D^{\alpha}(f(t)-f_m(t)) &=(g(t)-g_{m}(t))\\
&+\int_{0}^{t} k_{1}(t,s)f(s)-k_{1,m}(t,s)f_{m}(s) ds\\ 
&+ \int_{0}^{t} k_{2}(t,s)f(s)- k_{2,m}(t,s)f(s) dB(s) -H_{m}(t), t \in [0,1].\ (34)
\end{align*}
We rewrite equation (34) in the following form\\
\begin{equation}
e_m(t)=I^{\alpha}(\int_{0}^{t} k_1(s,t) e_{m}(s) ds)+ I^{\alpha}(\int_{0}^{t} k_2(s,t) e_{m}(s) dB(s))- H_m(t).
\end{equation} 
So, we have\\
\begin{equation}
e_{m}(t)=E_1+E_2+E_3.
\end{equation}
Now, we state some properties and theorem in the following.\\
Suppose that $f_{m}(t)$ is the best approximation of $f(t)$ which is approximated by shifted Legendre polynomials in $[a,b]$ so, we have\\
\textbf{Theorem 4}. Let $f_{m}(t)=\bar{f} \psi(t)$ where $\bar{f}=[f_0,\ldots,f_m]$ and 
\begin{equation*}
f_{m}=(2m+1) \int_{0}^{1} f(t) P_m(t) dt,
\end{equation*} 
then there exists a real number $\bar{z}$ such that
\begin{equation}
\|f(t)-f_m(t)\|_{2} \leq \bar{z} \frac{1}{(m+1)! 2^{2m+1}}
\end{equation}\\
where $\|f\| = \int_{0}^{1} |f(t)|^{2} dt$.\\
\textbf{Proof}. See [16]. \\
To approximate $E_1 and E_2$ by applying equation(37) we have\\
\begin{equation*}
\|E_1\|\leq \frac{1}{\Gamma(\alpha)} max_{s,t\in I}|k_1(s,t)| \|f(t)-f_m(t)\|,
\end{equation*}
so we have\\
\begin{equation}
\|E_1\|\leq \bar{z} \frac{1}{\Gamma(\alpha) (m+1)! 2^{2m+1}} max_{s,t\in I}|k_1(s,t)|.
\end{equation}
In the same way for $E_2$, we can write\\
\begin{equation}
\|E_2\|\leq \bar{z} \frac{1}{\Gamma(\alpha) (m+1)! 2^{2m+1}} max_{s,t\in I}|k_2(s,t)|
\end{equation}
and for $E_3$ we have\\
\begin{equation}
\|E_3\|=\frac{1}{\Gamma(\alpha)}\|H_m(t)\|.
\end{equation}\\
Therefore, by replacing equations(38-40) in equation (36), we obtain:
\begin{equation}
\|e_{m}\|^{2}\leq \frac{1}{Gamma(\alpha)}(\bar{z} \frac{max_{s,t \in I}(|k_1(s,t)|+|k-2(s,t)|)}{(m+1)!2^{2m+1}}).
\end{equation}

\section{Application and Numerical Examples}
In this section, some applications of stochastic differential equations, fractional stochastic differential equations and fractional Fokker-Planck equations are presented. Moreover, numerical examples are given to illustrate the high accuracy of the presented method.\\
In recent years, it has turned out many phenomena in signal processing and filtering, Population dynamics, economics, finance, fluid dynamics, microelectrons, structural mechanics, chemistry and medicine and etc can be successfully modelled by stochastic differential equations[11,9,12,1]. Ornstein and Uhlenbeck expressed a model to study the behaviour of gasses by the theory of stochastic differential equations. Ornstein and Uhlenbeck were interested in the velocity of an individual molecule of gas. Their equation is given by
\begin{equation*}
dX_{t}=-\alpha X_{t} dt+ \rho dB_{t}  with X_{0}=x_{0}
\end{equation*}
where $\alpha$ and $\rho$ are positive constants. In this equation $-\alpha X_{t} dt$ is negative when $X_{t}$ is larger than zero, and is positive when $X_{t}$ is smaller than zero[11].\\
Stochastic modelling has an important role in financial economics. For example, the Black-Scholes model [5] used abstract notations such as the It$\hat{o}$ calculus. The fractional Black-scholes market consist of a bank account or a bond and a stock also, by a fractional Volterra equations of the following form, we are able to govern the price $X_{t}$ of the stock 
\begin{equation*}
X_{t}=X_{0}+ \frac{1}{\Gamma(\beta)} \int_{0}^{t} \frac{\mu(s)}{(t-s)^{1-\beta}} ds+ \frac{1}{\Gamma(\frac{(\beta+1)}{2})} \int_{0}^{t} \frac{\chi(s)X_{s}}{(t-s)^{\frac{(1-\beta)}{2}} dW_{s}}
\end{equation*}
where ${W_{t}(\omega)=w(t), t>=0}$ is a Wiener process, $t\in[0,T]$, $0<\beta<1$ and the drift $\mu>=0$ and volatility $\chi>0$ are continuous functions on $[0,T]$.\\
 
\begin{figure}[!ht]
	\centering
	\includegraphics[width=55mm]{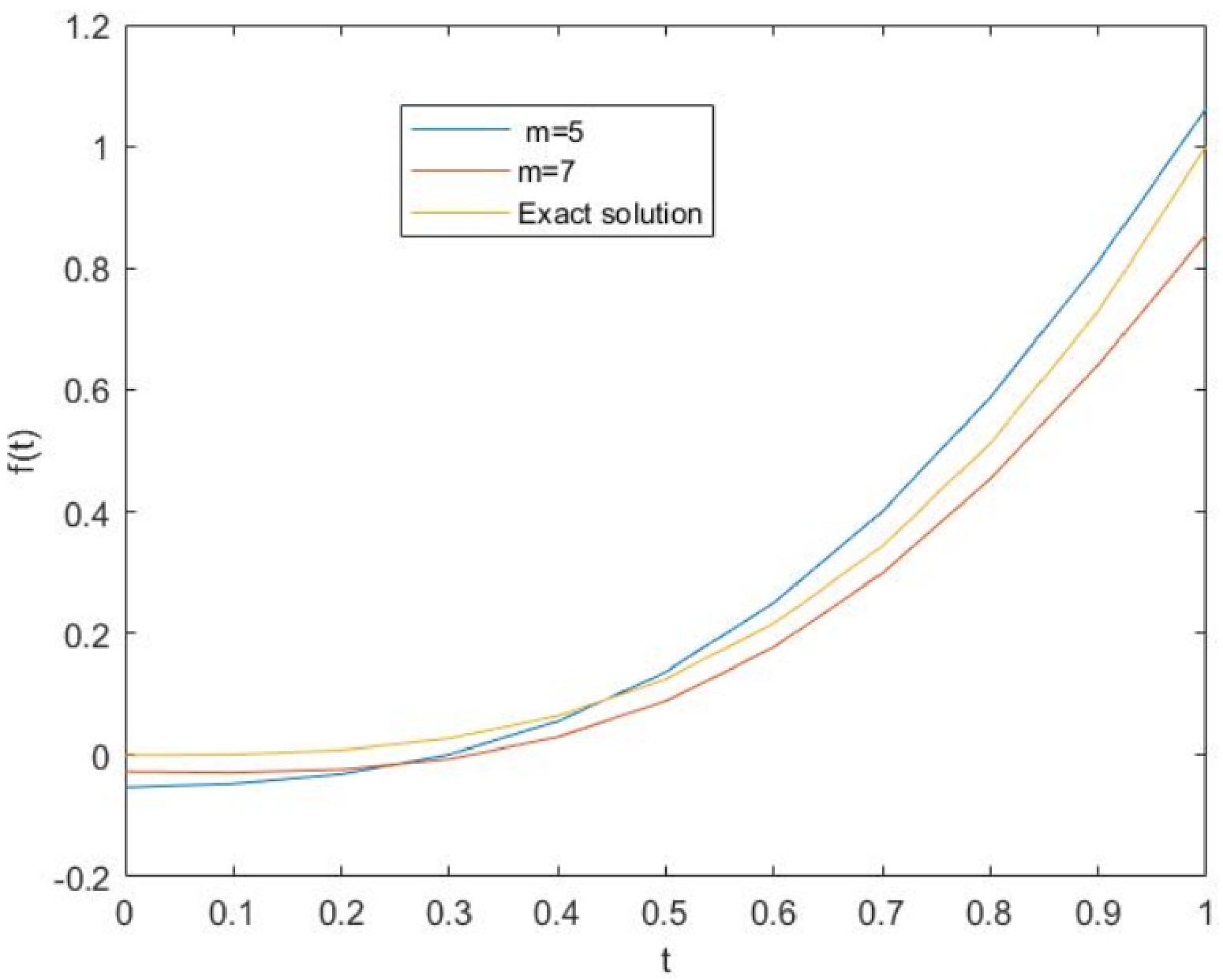}
	\label{fig:ex1}
	\caption{Solution with $\theta=0.75$, $\sigma=0$ is $f(t)=t^3$}
\end{figure}

\begin{figure}[!ht]
	\centering
	\includegraphics[width=55mm]{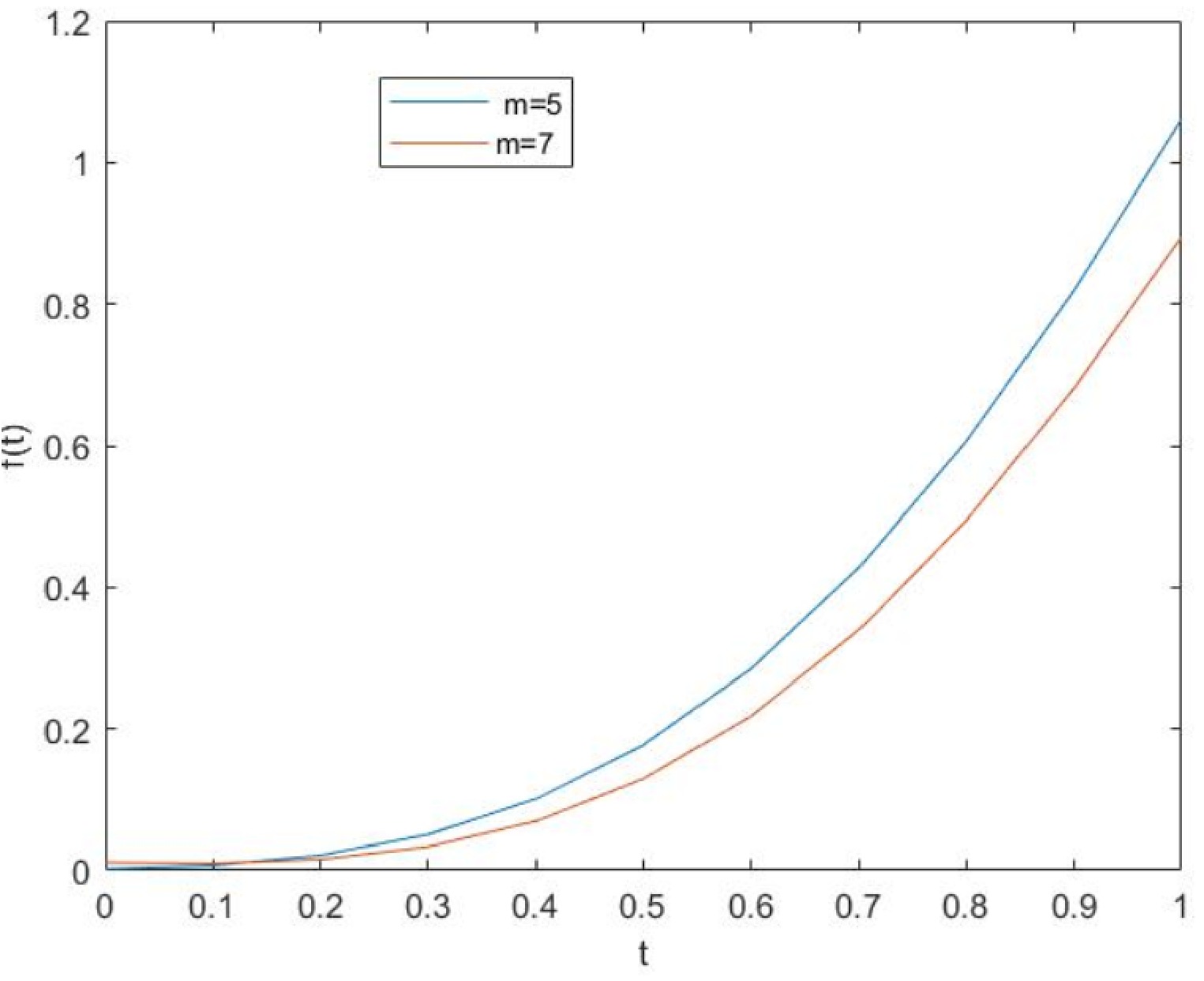}
	\label{fig:ex2v2}
	\caption{Solution with $\alpha=0.75, \sigma=1$ with $m=5,7$}
\end{figure}

\begin{example} 
Consider the following fractional stochastic integro-differential equation [38]
\begin{equation}
D^\alpha f(t)=-\frac{t^5 e^t}{5}+\frac{6t^{2.25}}{\Gamma(3.25)}+\int_{0}^{t} e^t s f(s) ds+ \sigma \int_{0}^{t} e^t s f(s) dB(s), t \in [0,1],
\end{equation}
where the exact solution with $\theta=0.75$, $\sigma=0$ is $f(t)=t^3$.\\
In Fig.~\ref{fig:ex1}, the results for $\alpha=0.75, \sigma=0$, by applying the technique described in section 4 with $m=5, 7$ are obtained and in Fig.~\ref{fig:ex2v2}, the results for $\alpha=0.75, \sigma=1$ with $m=5,7$ are obtained. Furthermore, to make a comparsion, you can see the approximate solution by the methods introduced in [38,9]. \qquad \textcolor{blue}{\Squaresteel}
\end{example}

\begin{figure}
	\centering
	\includegraphics[width=55mm]{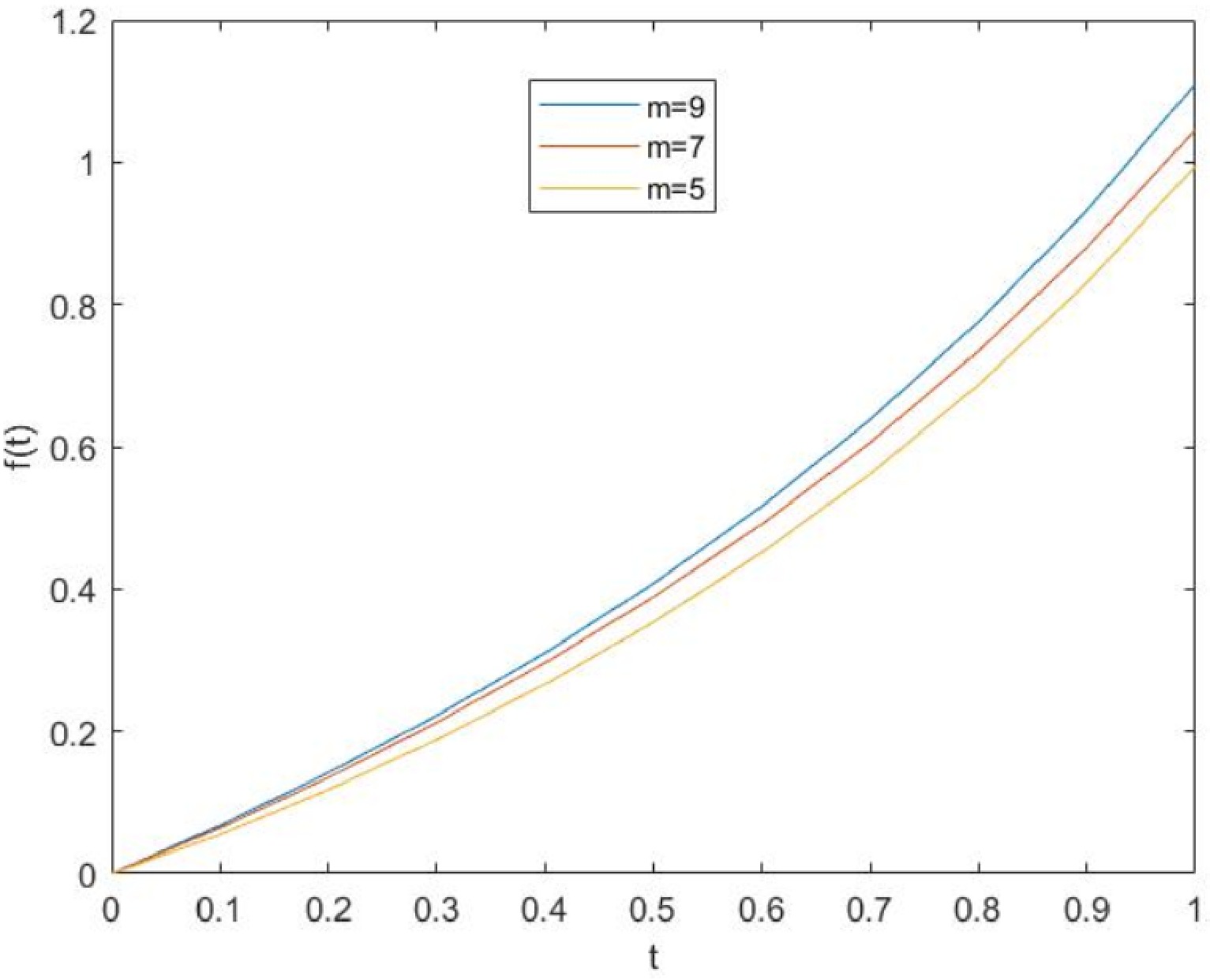}
	\label{fig:ex2}
	\caption{Solution with $\alpha=0.75, \sigma=1$ with $m=5,7,9$}
\end{figure}

\begin{example} 
Consider the following FSIDE [38]
\begin{equation}
D^\alpha f(t)=\frac{\Gamma(2)t^{1-\alpha}}{\Gamma(3-\alpha)} +\frac{t^3}{3}+ \int_{0}{t} sf(s) ds+ \int_{0}^{t} f(s) dB(s), t \in [0,1]
\end{equation}
with $f(0)=0$.\\
The results for example 2 are represented in Fig.~\ref{fig:ex2} with $m=5,7,9$ and $\alpha=0.75$.\qquad \textcolor{blue}{\Squaresteel}
\end{example}


\section*{Acknowledgements}
The research J.F. Peters has been supported by the Natural Sciences \&
Engineering Research Council of Canada (NSERC) discovery grant 185986 
and Instituto Nazionale di Alta Matematica (INdAM) Francesco Severi, Gruppo Nazionale per le Strutture Algebriche, Geometriche e Loro Applicazioni grant 9 920160 000362, n.prot U 2016/000036  and Scientific and Technological Research Council of Turkey (T\"{U}B\.{I}TAK) Scientific Human
Resources Development (BIDEB) under grant no: 2221-1059B211301223.
\section*{Declarations}
\section*{Ethics approval}
This article does not contain any studies with human participants or animals performed by any of the authors.\
\section*{Funding}
No funding to declare.
\section*{Competing interests}
The authors declare that they have no competing interests.
\section*{Authors' contributions}
All authors contributed equally to the writing of this paper. All authors read and approved the final manuscript.
\section*{Acknowledgements}
The authors express their deep gratitude to the referees for their valuable suggestions and comments for improvement of the paper.
\section*{References}
\noindent
\text{[1]} M. Asgari, Block pulse approximation for fractional stochastic integro-differential equations. Communications in Numerical Analysis, 2014(2014); 1-7.\\
\noindent
\text{[2]} H.A. Bahrawy, M.A. Alghamdi, A shifted Jacobi Gauss collocation scheme for solving fractional neutral functional-differential equations. Advances in Mathematical Physics, 2014(2014); 1-8.\\
\noindent
\text{[3]} J. Bao, B. Bottcher, X. Mao, C. Yuan,  Convergence rate of numerical solution to SFDEs with jumps. J. Comput. Appl. Math, 236(2011); 119-131.\\
\noindent
\text{[4]} G. Baumann, F. Stenger, Fractional Fokker-Planck equation. MDPI, Mathematics(2017); 1-19.\\
\noindent
\text{[5]} F. Black, M. Scholes, The prices of options and corporate liabilities. J.Political Economy, 81(1973); 637-654.\\
\noindent
\text{[6]} K. Burrage, P. Burrage, High strong order explicit Runge-Kutta methods for stochastic ordinary differential equations. Appl. Numer. Math, 22(1996); 81-101.\\
\noindent
\text{[7]} A. Dascioglu, A Chebyshev polynomials approximation for high-order partial differential equations with complicated conditions. Num. Math. Part. Diff. Equ, 25(2008); 610-619.\\
\noindent
\text{[8]} E. L. Doha, An accurate double Chebyshev spectral approximation for Poisson's equation, 1989.\\
\noindent
\text{[9]} M. EL-Borai, K. EL-Nadi, M. Osama, H. Ahmed, Volterra equations with fractional stochastic integrals. Mathematical Problems in Engineering, 5(2004); 453-468.\\
\noindent
\text{[10]} J. Franklin, Difference methods for stochastic ordinary differential equations. Math. Comput, 19(1965); 552-561.\\
\noindent
\text{[11]} J. A. Freund, T. P$\ddots{o}$schel, Stochastic processes in physics, chemistry, and biology. Springer- Velag Berlin Heidelberg, 2010.\\
\noindent
\text{[12]} R. Friedrich, J. Peinke, C. Renner, How to quality determinstic and random influences on the stochastics of the foreign exchange market. Phys. Rev. Lette., 84(2000); 5224-34.\\
\text{[13]} M. Gachpazan, M. Erfanian, H.  Beiglo, Solving nonlinear Volterra integro-differential equations by using Legendre polynomial approximations. Iranian Journal of Numerical Analysis and Optimization, 4(2014); 73-83.
\noindent
\text{[14]} E. Hesameddini, M. Shahbazi, Two-dimensional shifted Legendre polynomials operational matrix method for solving the two-dimensional integral equations of fractional order. Appl. Math. Copmut, 322(2018); 40-54.\\
\noindent
\text{[15]} N. Hofman, E. Platen, Stability of weak numerical schemes for stochastic differntial equations. Comput. Math. Appl, 28(1994); 45-57.\\
\noindent
\text{[16]} N. Hofman, Stability of weak numerical schemes for stochastic differential equations. Math. Comput. Simul, 38(1995); 63-68.\\
\noindent
\text{[17]} R. Janssen, Discretization of the Wiener-process in difference methods for stochastic differential equations. Stoch. Process. Their Appl, 18(1984); 361-369.\\
\noindent
\text{[18]} M. M. Khader, Numerical solution of nonlinear multi-order fractional differential equations by implementation of the operational matrix of fractional derivative. Studi. Nounlinear. Sc, 2(2011); 5-12.\\
\noindent
\text{[19]} M. Khodabin, K. Maleknejad, M. Rostami, M. Nouri, Numerical solution of stochastic differential equations by second order Runge-Kutta methods. Math. Comput. Model, 53(2011); 1910-1920.\\
\noindent
\text{[20]} L. Kloeden, R. Pearson, The numerical solution of stochastic differential equations. J. Aust. Math. Soc. Ser. B. Appl. Math, 20 (1992);8-12.\\
\noindent
\text{[21]} P. Koleden, E. Platen, A survey of numerical methods for stochastic differential equations. Stoch. Hydrol. Hydraul, 3(1989); 155-178.\\
\noindent
\text{[22]} S. Kouhkani, H. Koppelaar, Numerical solution of fractional differential equations by a shifted Chebyshev computational matrix. CJCME, 1(2018); 5-14.\\
\noindent
\text{[23]} S. Kouhkani, H. Koppelaar, R. Pettersson, Numerical solution of fractional stochastic integro-differential equations by the operational Tau method. InternationaL Journal of Statistical Analysis, 1(2019); 1-9.\\
\noindent
\text{[24]} S. Kouhkani, H. Koppelaar, M. Abri, Numrical solution of fractional neutral functional-differential equations by the Operational Tau Method. CJCME, 1(2017); 5-19.\\
\noindent
\text{[25]} X. Li, X. Mao, G. YIN, Explicit numerical approximations for stochastic differential equations in finite and infinite horizons: truncation methods, convergence in pth moment ans stability. IMA J. Numer. Anal, 1(2018); 1-46.\\
\noindent
\text{[26]} K. Maleknejad, S. Sohrabi, Legendre polynomial solution of nonlinear Volterre-Fredholm integral equations. IUST International Journal of Engineering Science, 19(2008); 49-52.\\
\noindent
\text{[27]} X. Mao, Stochastic differential equations and applications. 2nd ed., Ellis Horwood Publishing, Chichester, UK, 2007.\\
\noindent
\text{[28]} F. Mirzaee, N. Samadyar, On the numerical solution of fractional stochastic integro-differential equations via meshless discrete collocation method based on radial basis functions. Eng. Anal. Bound. Elem, 100(2018); 246-255.\\
\noindent
\text{[29]} T. Misawa, Numerical integration of stochastic differential equations by compoition methods. Dyn. Syst. Differ. Geom, 1180(2000); 166-190.\\
\noindent
\text{[30]} S. Nemati, P. M. Lima, Y. Ordokhani, Numerical solution of a class of two- dimensional nonlinear Voltera integral equations using Legendre polynomials. J. Comput. Appl. Math, 242(2013); 53-69.\\
\noindent
\text{[31]} J. F. Peters, Planar Ribbons, Ribbon Complexes and their Proximities. Ribbon Nerves, Betti Numbers and Descriptive Proximity, Bull. Allahabad Math. Soc., 2000, 1-13 and arXiv.1911.09014, 2019.\\
\noindent
\text{[32]} E. Platen, An introduction to numerical methods for stochastic differential equations. ActaNumer, 8(1999); 197-246.\\
\noindent
\text{[33]} A. Saadatmandi, M. Dehghan, A new operational matrix for solving fractional-order differential equations. Comput. Math. Appl, 59(2010); 1326-1336.\\
\noindent
\text{[34]} A. Saadatmandi, M. Dehghan, Numerical solution of the higher-order linear Fredholm integro-differential-difference equation with variable coefficients. Computers and Mathematics with Applications, 59(2010); 2996-3004.\\
\noindent
\text{[35]} I. Shehili, A. Zerroug, Bivariate Legendre approximation. Int. J. Appl. Math. Research, 6(2017); 125-129.\\
\noindent
\text{[36]} A. P. Smirnov, A. B. Shmelev, E. Ya. Sheinin, Analysis of Fokker-Planck approach for forein exchange market statistics study. Physica A, 344(2004); 203-206.\\
\noindent
\text{[37]} E. H. Spanier, Algebraic topology. McGraw-Hill Book Co., New York-Toronto, Ont.,CA, 1966, {x}iv+528 pp.,MR0210112. \\
\noindent
\text{[38]} Z. Taheri, S. Javadi, E. Babolian, Numerical solution of stochastic fractional integro-diffferential equation by the spectral collocation method. J. Comput. Appl. Math, 321(2017); 336-347.\\
\noindent
\text{[39]} G. Yuksel, M. Sezer, A Chebyshev series approximation for linear second-order partial differential equation with complicated conditions. Gaz. Uni. J. Sci, 26(2013);1-9.\\
\noindent
\text{[40]} S. Zhao, H. Jin, Strong convergence of implicit numerical methods for nonlinear stochastic functional differential equations. J. Comput. Appl. Math. 324(2017); 241-257.\\
\noindent
\text{[41]} J. F. Peters, T. Vergili, Fixed point property of amenable planar vortexes, Appl. Gen. Topol. 22 (2021), no. 2, 385--397, MR4359776.\\
\noindent
\text{[41]} J. F. Peters, Amiable and almost amiable fixed sets. Extension of the Brouwer fixed point theorem., Glas. Mat. Ser. III, 56(76) (2021), no.1, 175--194, MR4339174. \\
\noindent
\textbf{[42]} H. Thabet, S. Subhash, J.F. Peters, Numerical solutions for space–time conformable nonlinear partial differential equations via a scientific machine learning technique. Neurocomputing 620 (2025), 129134, 11--19.\\
\noindent
\textbf{[43]} H. Thabet, S. Subhash, J.F. Peters, Advances in solving conformable nonlinear partial differential equations and new exact wave solutions for Oskolkov-type equations. Math. Methods Appl. Sci. 45 (2022), no. 5, 2658–2673, MR4395619.\\
\noindent
\textbf{[44]} H. Thabet, S. Subhash, J.F. Peters, Analytical solutions for nonlinear systems of conformable space-time fractional partial differential equations via generalized fractional differential transform. Vietnam J. Math. 47 (2019), no. 2, 487–507, MR3986774\\
\noindent
\textbf{[45]} S. Smorgoni, Fractional Brownian Motion and Hurst Parameter Estimation in Stochastic Volatility Models, thesis, Universita degli Studi di Padova (2017).\\
\noindent
\textbf{[46]} D. Lehmann, The Dynamics of the Hubbard Model Through Stochastic Calculus and Girsanov Transformation, Int. J. of Theoretical Physics 63 (2024), 138--203, MR4752343.


\end{document}